\documentclass[conference]{IEEEtran}
\IEEEoverridecommandlockouts
\usepackage[keeplastbox]{flushend} 
\usepackage{array}
\usepackage{cite}
\usepackage{amsmath,amssymb,amsfonts}
\usepackage{algorithmic}
\usepackage{algorithm}
\usepackage{graphicx}
\usepackage{textcomp}
\usepackage{xcolor}
\usepackage{subcaption}
\usepackage{epstopdf}

\def\BibTeX{{\rm B\kern-.05em{\sc i\kern-.025em b}\kern-.08em
    T\kern-.1667em\lower.7ex\hbox{E}\kern-.125emX}}
\begin{document}

\title{Wireless Charging Lane Deployment in Urban Areas Considering Traffic Light and Regional Energy Supply-Demand Balance}

\author{\IEEEauthorblockN{Tian Wang, Bo Yang, Cailian Chen, Xinping Guan}}

\maketitle

\begin{abstract}
In this paper, to optimize the Wireless Charging Lane (WCL) deployment in urban areas, we focus on installation cost reduction while achieving regional balance of energy supply and demand, as well as vehicle continuous operability issues. In order to explore the characteristics of energy demand in various regions of the city, we first analyze the daily driving trajectory of taxis in different regions and find that the daily energy demand fluctuates to different degrees in different regions. Then, we establish the WCL power supply model to obtain the wireless charging supply situation in line with the real urban traffic condition, which is the first work considering the influence of traffic lights on charging situation. To ensure minimum deployment cost and to coordinate the contradiction between regional energy supply-demand balance and overall supply-demand matching, we formulate optimization problems ensuring the charge-energy consumption ratio of vehicles. In addition, we rank the priority of WCL efficiency to reduce the complexity of solution and solve the Mixed Integer NonLinear Programming (MINLP) problem to determine deployment plan. Compared with the baseline, the proposed method in this paper has significantly improved the effect.
\end{abstract}


\section{Introduction}
In recent years, with the promotion of Electric Vehicles (EVs), shared EV is also becoming a popular trip mode. However, due to the limitations of battery capacity and charging facilities, EV drivers generally face range anxiety and inconvenient charging. Traditional plug-in charging stations own relatively limited space and cause much recharging downtime, so they are difficult to satisfy the frequent charging needs of future vehicles (e.g., shared EVs). Hence dynamic charging based on Wireless Power Transfer (WPT) technique, by which moving vehicles are able to recharge on the Wireless Charging Lanes (WCLs), has gained much attention to overcome these battery-related issues \cite{b8}\cite{b9}, and made great progress at application in city \cite{b10}. 

It however brings up a problem for government: how to plan the deployment of WCLs to meet the frequent charging demand of EVs while reducing installation costs. These issues are fundamental in charging infrastructure installation scheme. To satisfy the travel needs of EVs, it is significant to analyze the real charging demand. Different functional regions in urban areas reflect different traffic conditions\cite{b14}, which probably determine various daily charging demand features. \textit{This is a factor worth considering}. As for dynamic charging, different vehicle objects and different scenarios will impact the consideration of deployment plan. Urban areas are the most suitable places to set up WCLs, and Shared EVs are the most beneficial objects. In urban areas, WCLs are usually considered to be located at intersections for covering most vehicles\cite{b1} \textendash\cite{b13}. Hence, Traffic Light (TL) is a significant factor that affects vehicle charging conditions, while \textit{there have been no previous works that consider TL in dynamic charging model}.  

There has been some researches focus on the optimal placement of charging stations for Plug-in EVs to improve the driving experience \cite{b11}\cite{b12}. Generally, the point of views they consider from can be divided into three aspects: financial cost (build up cost\cite{b15} and substation power loss cost\cite{b16}), the minimum average distance between charging stations and vehicles\cite{b11}\cite{b15}, and the reduction in average charging time\cite{b12} \cite{b18}. However, the plug-in charging station deployment schemes cannot be directly applied to the WCL deployment given the following difference: charging approach of the former is not closely related to the traffic condition, while the one of the latter is deeply impact by the behavior of vehicle (e.g., speed) and state of traffic (e.g., traffic flow and TL situation).

Although a few researches have also paid attentions on the optimal scheme of wireless charging infrastructures allocation, there are still some aspects that have not been studied. Team of Korea Advanced Institute of Science and Technology has focused on their \textit{OnLine Electric Vehicle} (OLEV) for several years: the economic benefit of dynamic charging EV is analyzed in \cite{b10}, and system architecture and mathematical models of electric transit bus system model is built to optimally allocate power tracks \cite{b5}, while determining the battery size of EV in multiple route environment \cite{b8}. But their researches are only based on OLEV bus system, lacking in applicability for other vehicle types. In \cite{b18}, a method for optimal placement of WPT system for autonomous driving vehicles is proposed from a viewpoint of optimal control problem. In \cite{b2}, multiple sources of vehicles in urban areas are considered, and by `categorization-clustering-extraction-deployment' steps, the integer programming allocation problem are formulated and solved. Nevertheless, all of them ignore impacts of TL on charging model, which is inevitable in real urban areas. And the balance of supply and demand have not been studied, which is one of the main reasons the government builds dynamic charging infrastructures.

Hence in this paper, we mainly aim at the reasonable power supply model that in line with actual traffic condition in urban areas, as well as the supply-demand balance for the characteristics of different functional regions. The contributions of this paper are summarized as follows. 
\begin{itemize}
	\item We find the fluctuation features of daily charging demand in the whole urban areas and different subregions, and formulate bilevel utility functions to resolve the contradictions between optimum of subregions and global optimum.
	\item We derive an hourly based WCL power supply model in line with the actual urban traffic according to the movement and queuing of vehicles in the TL cycle, and verify the suitability of this model with the result of supply and demand balance. 
	\item We transform the original Mixed Integer NonLinear Programming (MINLP) problem into a standard solvable Mixed Integer Linear Programming (MILP) problem, and solve it based on Generalized Benders Decomposition (GBD). Comparing proposed method with the baseline, we get a conclusion that lanes with maximum traffic flow are not entirely the optimal placement choices in urban areas, and industrial subregions are unsuitable for WCL allocation.
\end{itemize} 

The remainder of this paper is organized as follows. Section ~\ref{sct2} analyzes the fluctuation features of charging demand in urban subregions. Section~\ref{sct3} shows the WCL power supply model and problem formulation. The algorithm for resolving the MINLP issue is provided in Section~\ref{sct4}. Simulations to validate our results are given in Section~\ref{sct5}. And Section~\ref{sct6} provides the conclusion.
\section{Metropolitan-Scale Trajectory Dataset Analysis}
\label{sct2}
The purpose of building charging facilities is to meet the charging demand of shared EVs in the city, therefore, we first study charging demand of vehicles. While in general, characteristics of different functional regions of the city might vary, we study the daily charging demands of vehicles in different functional regions from real datasets.

\begin{figure}[t]
	\setlength{\belowcaptionskip}{-0.5cm}   
	\centerline{\includegraphics[width=0.9\linewidth]{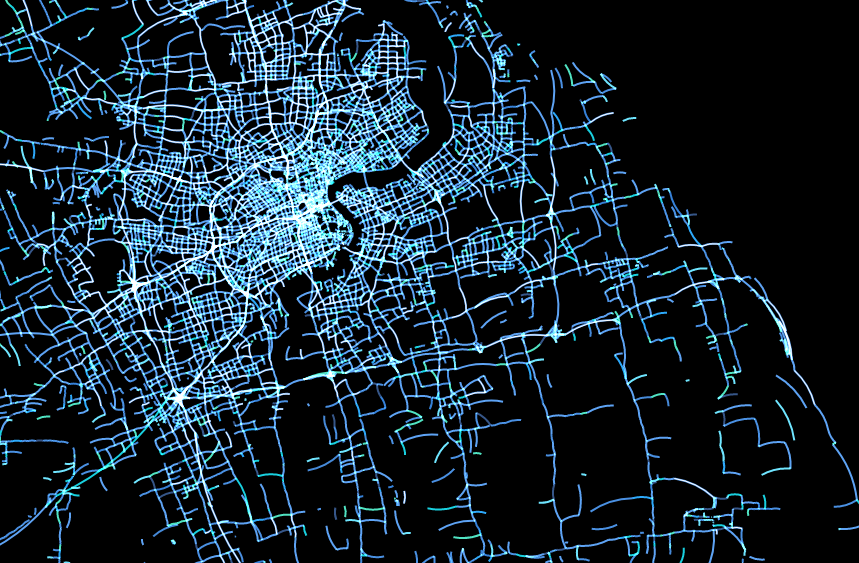}}
	\caption{Traffic heat map.}
	\label{heatmap}
\end{figure}

Our dataset, which is from Shanghai Traffic Information Center, records the GPS trajectory and speed data of about 13,000 taxis in Shanghai for one month (April 1$\sim$30, 2015), with a 10-second recording period. As the driving trajectory of a vehicle directly reflects its energy consumption, and the daily trip mode of taxis is similar to that of shared EVs, their trajectory data can well reflect the driving situation of the latter. Based on the above datasets, we drew a city traffic heat map to observe the traffic conditions in urban areas. As Fig.~\ref{heatmap} shows, the brighter roads represent that taxis on which have lower average speed, in other words, they are more congested.

\begin{figure}[t]
	\setlength{\belowcaptionskip}{-0.5cm}  
	\centerline{\includegraphics[width=0.9\linewidth]{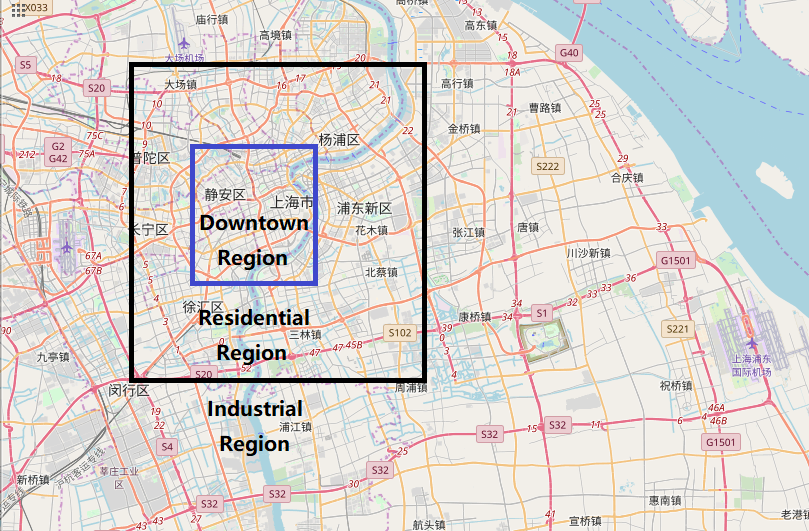}}
	\caption{Function regions.}
	\label{regions}
\end{figure}

\begin{figure*}[htbp]
	\begin{minipage}[t]{0.245\linewidth}
		\centerline{\includegraphics[width=1\textwidth]{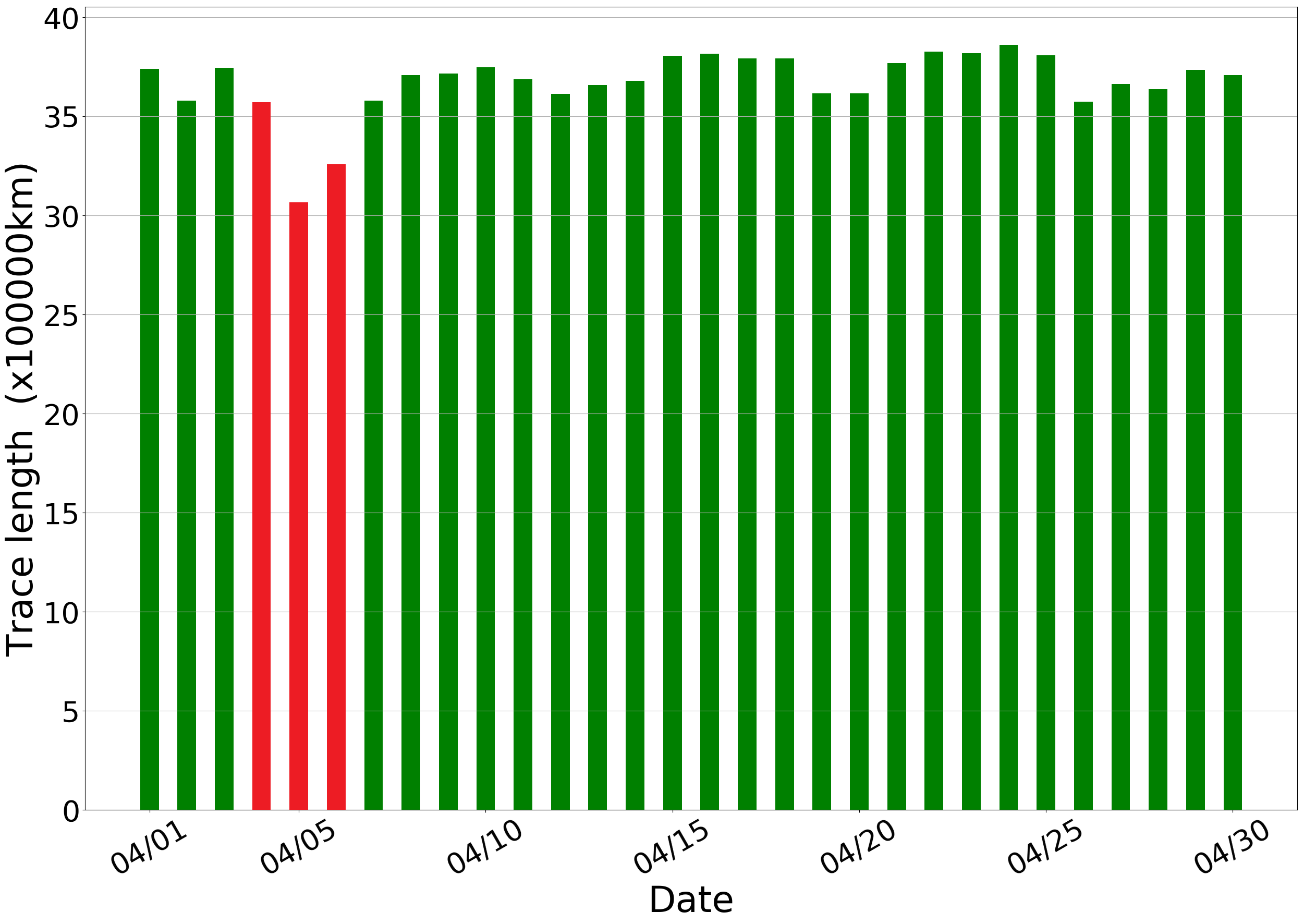}}
		\setlength{\abovecaptionskip}{-1pt}
		\subcaption{Total urban areas.}
		\label{all}
	\end{minipage}
	\begin{minipage}[t]{0.245\linewidth}
		\centerline{\includegraphics[width=1\textwidth]{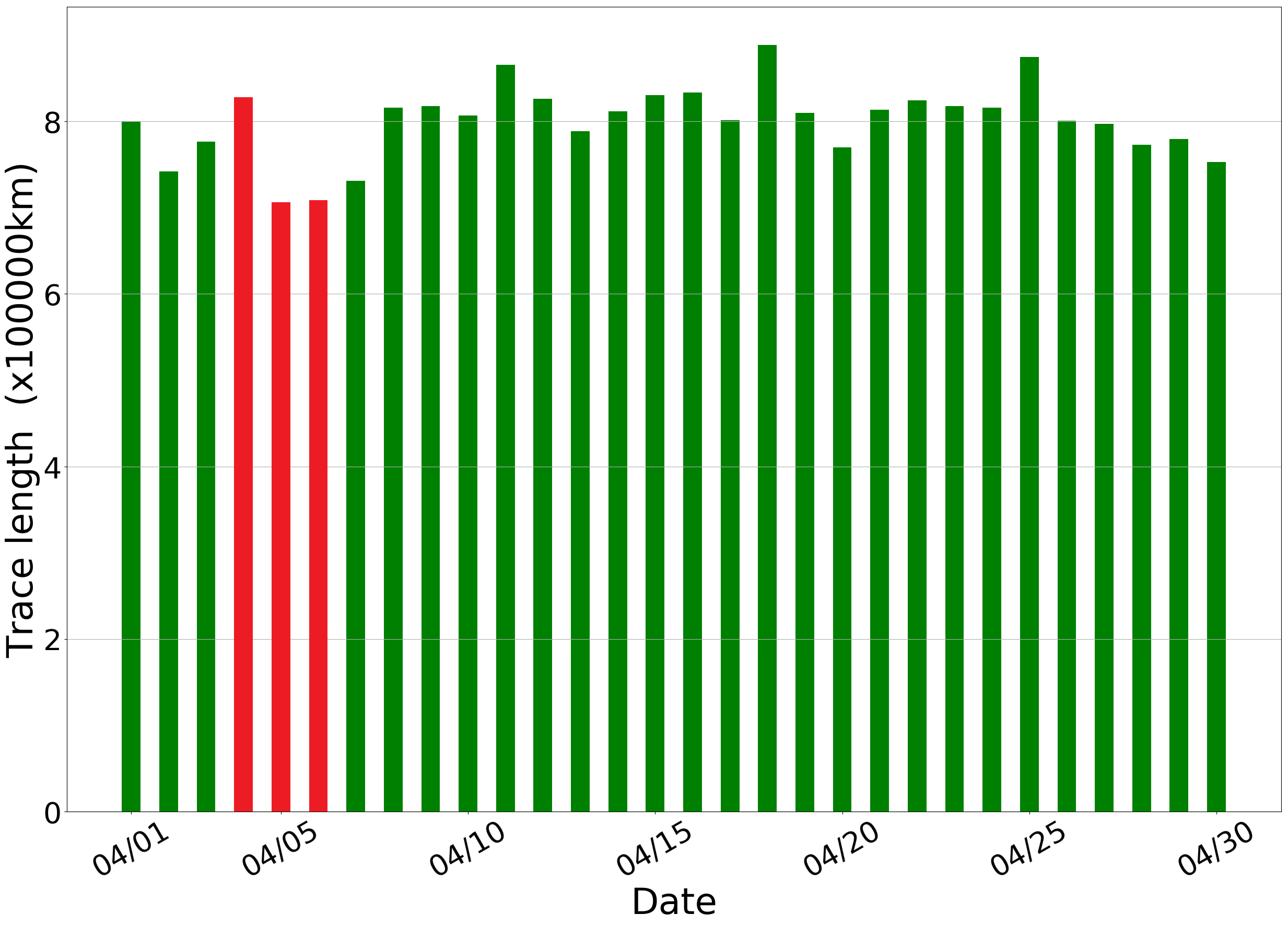}}
		\setlength{\abovecaptionskip}{-1pt}
		\subcaption{Downtown region.}
		\label{downtown}
	\end{minipage}
	\begin{minipage}[t]{0.245\linewidth}
		\centerline{\includegraphics[width=1\textwidth]{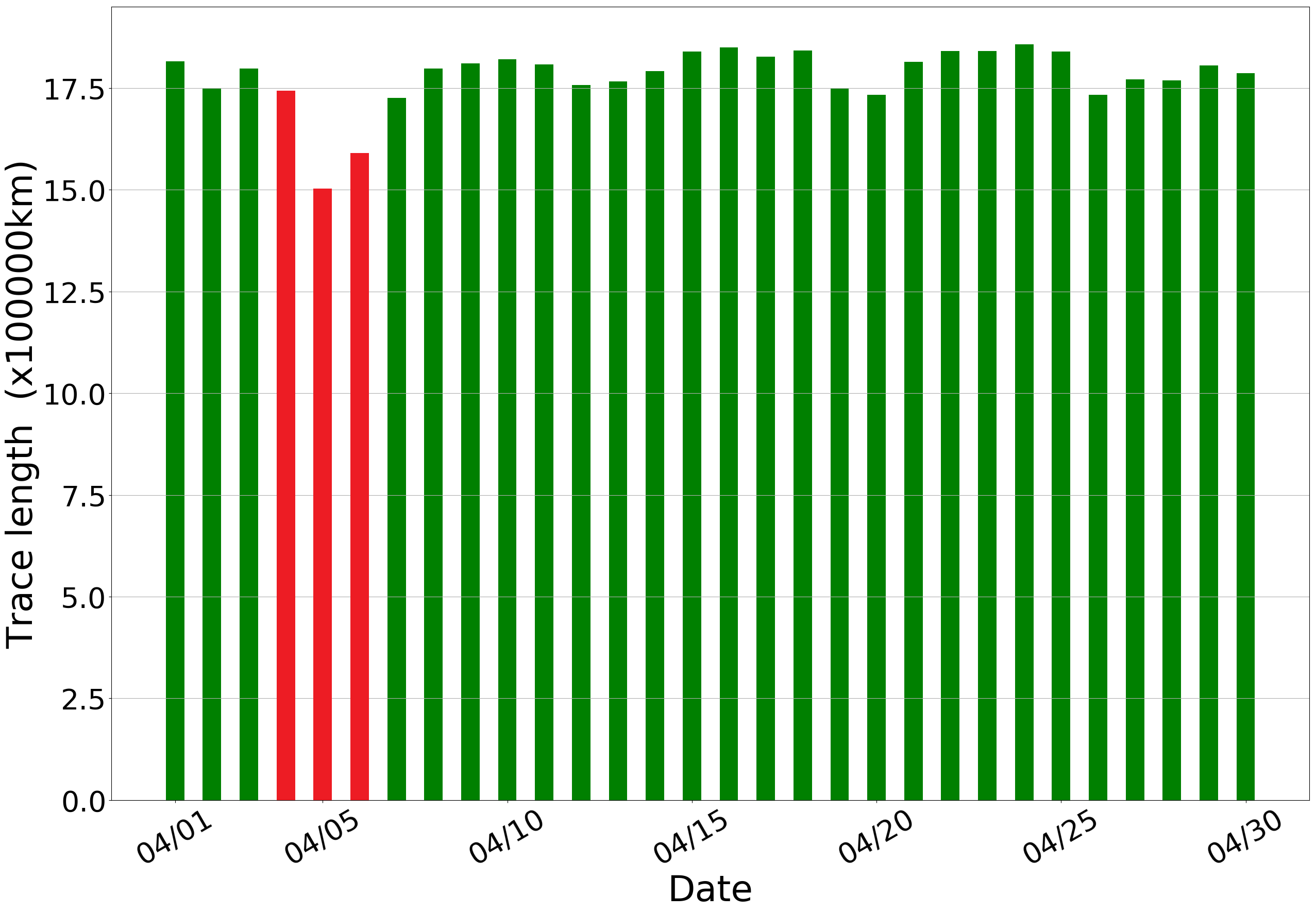}}
		\setlength{\abovecaptionskip}{-1pt}
		\subcaption{Residential region.}
		\label{residensial}
	\end{minipage}
	\begin{minipage}[t]{0.245\linewidth}
		\centerline{\includegraphics[width=1\textwidth]{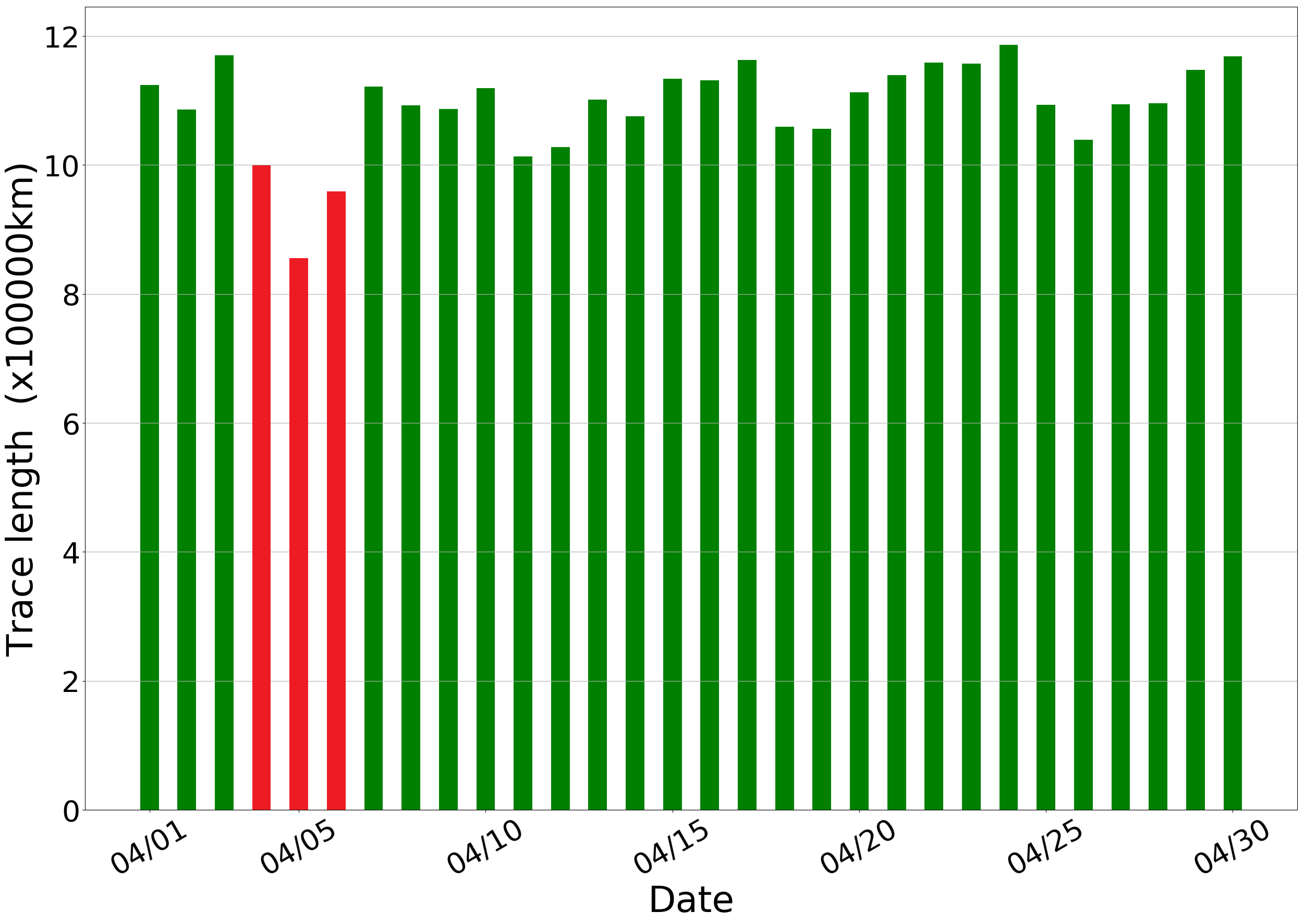}}
		\setlength{\abovecaptionskip}{-1pt}
		\subcaption{Industrial region.}
		\label{industrial}
	\end{minipage}
	\setlength{\belowcaptionskip}{-0.3cm}  
	\caption{Daily trajectory in regions.}
	\label{trajectory}
\end{figure*}

According to the traffic heat map and landmarks in city, the main urban areas are simply divided into three functional subregions: \textit{Downtown region}, \textit{Residential region} and \textit{Industrial region} (see Fig.~\ref{regions}). As indicated in Fig.~\ref{trajectory}, we calculated the total length of the daily trajectories in each functional subregion as well as in the whole urban areas for 30 days. The red bars means April 4$\sim$6 are \emph{Tomb Sweeping Day}, so we do not consider their impact. 
It is clearly that fluctuations of daily trajectory length in the whole urban areas and residential region are small, while the trajectory length in downtown region and industrial region vary obviously with seven days a week. Since the EV battery discharging quantity and driving distance can be considered as proportional in a subregion \cite{b2}\cite{b17}, and consumption is right demand, the charging demand in each subregion can be modeled as \eqref{eq23}.
\begin{equation}
N_i=q_i\cdot tr_i,\quad \forall i\in 1,\cdots, I,
\label{eq23}
\end{equation}
where the urban areas are divided into $I$ functional subregions, and $q_i$ represents the rate of EV battery discharge per unit of distance and $tr_i$ denotes the total trace length of EVs through $i$th subregion.

Hence it can be figured out from Fig.~\ref{trajectory} that the daily total charging demand of vehicles throughout all regions are essentially the same, and the energy supply of all the WPT lanes in the city should match the average charging demand, i.e., supply are supposed to equals to demand. However, the daily total charging requirements of vehicles in subregions have different fluctuation features, the average value cannot represents the real energy demand so well and the charging supply of WPT lanes in each region should approach to the average charging demand, i.e., they're not necessarily equal to each other. 

We assign normalized Mean-Square Error (MSE) to measure the fluctuation of the daily trajectory length, and results are shown in TABLE~\ref{tab1}. 
\renewcommand\arraystretch{1.5}
\begin{table}[htbp]
	\setlength{\abovecaptionskip}{-0.1cm} 
	\caption{Normalized MSE of Regions}
	\begin{center}
		\begin{tabular}{|c|c|c|c|c|}
			\hline
			\textbf{Regions} & \textbf{\textit{Total}}& \textbf{\textit{Downtown}}& \textbf{\textit{Rsidential}}& \textbf{\textit{Industrial}} \\
			\hline
			\textbf{MSE}& 0.0225& 0.0444& 0.0214& 0.0403\\
			\hline
		\end{tabular}
		\label{tab1}
	\end{center}
\end{table}
\section{Preliminaries And Problem Formulation}
\label{sct3}
In this section, we first describe the waiting queue length model, and then present the power supply capability model. In addition, we rank the candidate lanes, and finally formulate optimization problems to determinate the charging lane location. 

\subsection{Waiting Queue Length at Intersection}
The locations of charging lanes should consider two main factors: vehicle coverage and charging amount. The amount of charging is mainly determined by the charging duration when the charging power is constant. In an urban environment, it is obvious that intersections cover the traffic maximally, and cars move slowly and often stop at TLs\cite{b1}\cite{b2}, so we choose traffic intersections as candidate locations. 

Charging lane length determines how many cars can be charging simultaneously, so it needs to be optimized by considering waiting queue length (WQL). According to different traffic conditions, the average WQL at the intersection will be different. On one hand, the average speed of vehicles $v_{i,k}$ can partly reflect the traffic congestion. On the other hand, the average WQL fluctuates up and down near the mean number of passed vehicles in a single TL period. The worse traffic condition is, the greater WQL is, and vice versa. Hence we infer that such trend is approximately exponential, and the relationship between WQL and average speed can be modeled as:
\begin{equation}
	\setlength{\abovedisplayskip}{3pt}
	\setlength{\belowdisplayskip}{3pt}
	L^{w}_{i,k}=\dfrac{f_{i,k}}{n_{i,k}} e^{\varphi (\dfrac{1}{v_{i,k}}-\dfrac{1}{v_{equal}})}\label{eq1},\quad \forall k\in 1,\cdots,m_i
\end{equation}
Note that $f_{i,k}$ is traffic flow (the number of passing vehicles in a time slot, i.e., 1 hour), and $f_{i,k}/n_{i,k}$ is mean number of passed vehicles in a TL period (red light duration $T^{r}_{i,k}$ and green light duration $T^{g}_{i,k}$), which is called standard value. $\varphi$ is a constant parameter in a relative stable scenario. $m_i$ is the total number of deployed WCLs in the $i$th subregion.  Periods of TLs are assumed unaltered in a time slot, so we can calculate period numbers in a time slot $\tau$:
\begin{equation}
	\setlength{\abovedisplayskip}{3pt}
	\setlength{\belowdisplayskip}{3pt}
	n_{i,k}=\dfrac{\tau}{T^{r}_{i,k}+T^{g}_{i,k}}\label{eq2}
\end{equation}

\subsection{Power Supply Capability of WCL}
The power supply capability of WCL determines the total electric energy that facilities can provide for vehicles driving on the lane in a day. The charging situation is deeply affected by the traffic conditions, including the traffic flow, speed of vehicles, inter-vehicle distance and TL period. In above items, traffic flow and the average vehicular speed can be calculated by sensors, then the average inter-vehicle distance can be deduced from the relationship of traffic flow, traffic density and vehicular speed\cite{b4}.
\begin{equation}
	\bar{d_{1}}=\dfrac{1}{K}=\dfrac{v}{f},
\end{equation}
where $K$ is the traffic density, and $\bar{d_{1}}$ is the average headway of moving vehicles, which is constant in a time slot. 

To calculate the total power supply, we separately discuss the total charging amount of vehicles that are charged statically ever and ones only charge in transit. 

\subsubsection{Charge Statically Ever} 
When the TL turns red, stopped cars can utilize the standby time for recharging. Therefore we only consider vehicles that will stop in this duration. For convenience of analysis, WQLs in each red TL period are assumed as average WQL, and queuing cars pass through the intersection with a minimum safe speed $v_{min}$ when the green TL is on\cite{b3}. As \eqref{eq5} shows, to the $r$th EV in the queue, the charging duration consists of two intervals: interval before departure (first two terms) and departure interval.
\begin{equation}
	T=T^{r}_{i,k}+\dfrac{L_{i,k}-r\bar{d_1}}{v_{i,k}}+\dfrac{r\bar{d_2}}{v_{min}/2}
	\label{eq5}
\end{equation}
where $\bar{d_2}$ is the average headway of stopped vehicles.

We infer that the power of wireless charging facilities are constant, so charge quantity can be calculated. By accumulating the total charging amount of vehicles in each time slot, we can get the total power supply $C^1$ of a WCL. However, whether the length of WCL is longer than WQL leads to different power supply model.

\paragraph{WCL shorter than WQL ($L_{i,k}<L^w_{i,k}$)} 
\begin{equation}
	\begin{split}
		C^1_{i,k}&(L_{i,k})=pn_{i,k}\sum_{r=1}^{r_1}T\\
		&=pn_{i,k}\Big[r_2\Big(T^r_{i,k}+\dfrac{L_{i,k}}{v_{i,k}}\Big)+\sum_{r=1}^{r_2}\Big(\dfrac{\bar{d_2}}{v_{min}/2}-\dfrac{\bar{d_1}}{v_{i,k}}\Big)r\Big]\\
		&=pn_{i,k}\Big[\dfrac{L_{i,k}}{\bar{d_2}}\Big(T^r_{i,k}+\dfrac{L_{i,k}}{v_{i,k}}\Big)+\dfrac{1}{2}\Big(1+\dfrac{L_{i,k}}{\bar{d_2}}\Big)\dfrac{L_{i,k}}{\bar{d_2}}\xi\Big]\\
		&=pn_{i,k}\Big[\Big(\dfrac{\xi}{2\bar{d_2}^2}+\dfrac{1}{v_{i,k}\bar{d_2}}\Big)L^2_{i,k}+\Big(\dfrac{T^r_{i,k}+\xi/2}{\bar{d_2}}\Big)L_{i,k}\Big]
	\end{split}
	\label{eq7}
\end{equation}
\paragraph{WCL longer than WQL ($L_{i,k}\ge L^w_{i,k}$)}
\begin{equation}
	\begin{split}
		C^1_{i,k}&(L_{i,k})=pn_{i,k}\sum_{r=1}^{s}T\\
		&=pn\Big[s\Big(T^r_{i,k}+\dfrac{L_{i,k}}{v_{i,k}}\Big)+\sum_{r=1}^{s}\Big(\dfrac{\bar{d_2}}{v_{min}/2}-\dfrac{\bar{d_1}}{v_{i,k}}\Big)r\Big]\\
		&=pn\Big[s\Big(T^r_{i,k}+\dfrac{L_{i,k}}{v_{i,k}}\Big)+\dfrac{(1+s)s}{2}\xi\Big]\\
		&=\dfrac{pn_{i,k}s}{v_{i,k}}L_{i,k}+pn_{i,k}s\Big(T^r_{i,k}+\dfrac{(1+s)}{2}\xi\Big)
	\end{split}
	\label{eq8}
\end{equation}
where $\xi$ equals $\dfrac{\bar{d_2}}{v_{min}/2}-\dfrac{\bar{d_1}}{v_{i,k}}$, and $p$ is the constant charging power. Note that $r_1$ is the number of vehicles that WCL can contain and equals  $L_{i,k}/\bar{d_2}\ rounded\ down$, and $s$ is the number of the waiting vehicles and equals $L_{i,k}/\bar{d_1}\ rounded\ down$. Combining \eqref{eq7} and\eqref{eq8}, the model can be simplified as \eqref{eq9}.
\begin{equation}
	C^1_{i,k}(L_{i,k})=\left\{
	\begin{array}{rl}
		a_{i,k}L_{i,k}+b_{i,k}, &{\rm if}~L_{i,k}\ge L^w_{i,k}\\
		c_{i,k}L^2_{i,k}+d_{i,k}L_{i,k}, &{\rm otherwise}
	\end{array}
	\right.
	\label{eq9}
\end{equation}

To maximize charging efficiency, the length of WCL are supposed to cover the waiting queue, so that all the stopped vehicles are able to maximally utilize their standby time. Therefore, we only consider the situation when $L_{i,k}\ge L^w_{i,k}$.

\subsubsection{Only Charge in Transit}
When the traffic is not too heavy, some vehicles do not have to stop for the red TL and only recharge in transit; while in a terrible traffic condition, the road head cannot digest the waiting queue in a TL period. By comparing the total traffic flow and the number of queuing cars in a time slot, the model of $C^2_{i,k}$ is shown as \eqref{eq10}.
\begin{equation}	
	C^2_{i,k}(L_{i,k})=\left\{
	\begin{array}{rcl}
		p(f-n_{i,k}s)\dfrac{L_{i,k}}{v_{i,k}}, &{\rm if}~f\ge n_{i,k}s\\
		0, &{\rm otherwise}
	\end{array}\right.
	\label{eq10}	
\end{equation}

\subsubsection{Total Power Supply in a Region}
Combining \eqref{eq9} and \eqref{eq10}, we finally get the power supply capacity of a WCL as \eqref{eq11}. If the WCL locations are fixed, the total energy provided by WCLs in a region in unit time is obtained as \eqref{eq12} shows.
\begin{equation}
	C_{i,k}=C^1_{i,k}+C^2_{i,k}=a'_{i,k}L_{i,k}+b_{i,k}
	\label{eq11}
\end{equation}
\begin{equation}
	C_i=\sum_{k=1}^{m_i}C_{i,k}\cdot x_{i,k}=(L_iA_i+B_i)X^\top_i
	\label{eq12}
\end{equation}
where $L_i=\{L_{i,1},\cdots,L_{i,m_i}\}$ are the lengths of WCLs and $X_i=\{x_{i,1},\cdots,x_{i,m_i}\}$ are 0/1 variables representing whether the lanes are paved with WCL. $A_i$ is a diagonal matrix consisting of coefficients of $L_i$ and $B_i$ is a row vector.
\subsection{Charging Lane Location Determination}
\subsubsection{Sub-Utility Function}
In each subregion, we hope that the deployment cost, which is proportional to the WCL length\cite{b2}, can be minimized. At the same time, the gap between supply and demand is desired to be as small as possible. We use a tolerance coefficient $\gamma$ to control the gap term. Since it is a waste of resources to lay a WCL too long, we limit the lengths of WCLs in a range of $\varepsilon L^w_{i,k}$.
\begin{align}
	&\min_{L_{i,k}, x_{i,k}}\quad u_i=\sum_{k=1}^{m_i} \rho\cdot L_{i,k}\cdot x_{i,k}+\gamma_i|C_i-N_i|\label{eq13}\\	
		& \begin{array}{r@{\quad}r@{~}l}
		s.t.&L_{min}&\le L_{i,k}\le L_{max}\\
		&L^w_{i,k}&\le L_{i,k}\le \varepsilon L^w_{i,k}\\
		&x_{i,k}&\in\{0,1\}\\
		&\forall i&\in 1,\cdots,I,~k\in 1,\cdots,m_i
		\label{eq14}
	\end{array}
\end{align}	
where $\rho$ is the unit price per meter of WCL.
\subsubsection{Main-Utility Function}
We hope every subregion can maximize its sub-utility function. Therefore, we define a utility coefficient $uc$ to normalize the utility satisfaction, as \eqref{eq15} shows. The main problem is set as P1.
\begin{equation}
	uc_i=\dfrac{u_i(L_i,X_i)}{u^{min}_i}
	\label{eq15}
\end{equation}
\begin{equation}
	{\rm P1:} \quad\min_{L_i,X_i}\quad max\{uc_1,uc_2,\cdots,uc_I\}
\end{equation}

In case that a vehicle cannot reach an alternative WCL in inter-region trips with its state of charge (SOC), we consider the recharge-consumption ratio $R(i)$, which represents the ratio of average rechargeable energy $c_j$ to energy consumption $q_{ij}d_{ij}$ when a vehicle arrives at the next WCL $j$ (e.g., the $n$th one in the $m$th region) from intersection $i$, to ensure operability of EVs (maintaining enough residual power for continuously moving). To cover the possibility of passing the next WCL equipped intersections $in_j\in\widetilde{IN},$ where $\widetilde{IN}=\{in_1,in_2,\cdots,|\widetilde{IN}|\}$, we consider the weighted average ratio based on visiting probability $\hat{f}(in_j)$ as \eqref{eq18}. Assign $IN$ as all intersections in road network, the average recharge-consumption ratio at an intersection $in$ is:
\begin{equation}
R(i)=\dfrac{c_j}{q_{ij}d_{ij}}=\dfrac{C_{m,n}/f_{m,n}}{q_{ij}d_{ij}}
\end{equation}
\begin{equation}
\hat{f}(in_j)=\dfrac{f_j}{\sum_{k=1}^{|\widetilde{IN}|}f_k}
\label{eq18}
\end{equation}
\begin{equation}
\bar{R}(in)=\sum_{j=1}^{|\widetilde{IN}|}\hat{f}(j)R(in)x_j\quad ,\forall in\in IN
\end{equation}

By guaranteeing this ratio higher than a threshold $\alpha$, vehicles' energy supply is guaranteed within a certain distance. For the convenience of solving the main optimization problem, we convert P1 to P2:
\begin{align}
	{\rm P2:} \quad&\min_{L_i, X_i}\qquad\qquad\qquad\,V\\
	&\begin{array}{r@{\quad}r@{~}l}
		s.t.&\sum_{i=1}^{I}C_i&=\sum_{i=1}^{I}N_i=N\\
		&uc_i(L_i,X_i)&\le V\\
		&\bar{R}(in)&\ge\alpha,~\forall in\in|IN|\\
		&\theta(L_i,X_i)&\le 0\\
	\end{array}
\end{align}
where the first constraint represents that the total power supply in the whole urban area matches the total energy demand, and $\theta(L_i,X_i)\le 0$ is the constraint of subproblems.

\subsubsection{Rank of Candidate Locations}
To simplify calculation, we take the weighted traffic flow and power supply of the lane as the reference to design the priority rank of lanes as \eqref{eq22}, and sort all the lanes based on the rank. Finally the first $\eta$\% lanes are selected as the candidate locations. 
\begin{equation}
	E_{i,k}=\dfrac{f_{i,k}}{f^{max}_{i,k}}+\beta\dfrac{C_{i,k}}{C^{max}_{i,k}}
	\label{eq22}
	\setlength\abovedisplayskip{0.5cm}
	\setlength\belowdisplayskip{0.8cm}
\end{equation}

\section{Algorithm}
\label{sct4}
\subsection{Problem Transformation}
In \eqref{eq13}, the sub-utility function is a MINLP problem. However, current solvers cannot resolve MINLP problem with multiplication of continuous and integer variables like the first term. To decouple $L_i$ and $X_i$, we transform \eqref{eq12} into a solvable MILP formulation:
\begin{equation}
	C'_i=L_iA_iE+B_iX^\top_i,
\end{equation}
where E is a column vector whose elements are 1.

if $C'_i\ge N_i$:
\begin{align}
	&\min_{L_i, X_i}\quad u'_{i+}=L_i(\gamma_iA_iE+\rho E)+\gamma_iB_iX^\top_i-\gamma_iN_i\label{eq25}\\
	& \begin{array}{r@{\quad}r@{~}l}
		s.t.&\qquad L^w_{i,k}\cdot x_{i,k}&\le L_{i,k}\le \varepsilon L^w_{i,k}\cdot x_{i,k}\\
		&\qquad L_i(\gamma_iA_i&E+B_iX^\top_i)\ge N_i\\
	&&\quad \cdots
	\label{eq26}
	\end{array}
\end{align}

if $C'_i<N_i$, \eqref{eq25} is replaced with \eqref{eq27}:
\begin{equation}
	\min_{L_i, X_i}\quad u'_{i-}=L_i(\rho E-\gamma_iA_iE)-\gamma_iB_iX^\top_i+\gamma_iN_i\label{eq27}
\end{equation}	
The minimum value of \eqref{eq25} and \eqref{eq27} is optimal result.

In main-utility function, the original formation is also unsolvable. We introduce new variables $Q_i$ to enable utility function to be composed of continuous and integer variables, so that it can be decomposed later. 
\begin{equation}
	\begin{split}
		Q_i&=|\gamma_iL_iA_iE+\gamma_iB_iX^\top_i-\gamma_iN_i|+\rho L_iE-u^{min}_iV\\
		&=u'_i(L_i,X_i)-u^{min}_iV
	\end{split}
\end{equation}

The new formulation shows as \eqref{eq29} and \eqref{eq30}:
\begin{align}
	{\rm P2':} \quad&\min_{L_i, X_i}\quad \sum_{i=1}^{I}u'_i(L_i,X_i)-Q_i\label{eq29}\\
	&\begin{array}{r@{\quad}r@{~}l}
		s.t.&\sum_{i=1}^{I}C'_i&=\sum_{i=1}^{I}N_i=N\\
		&u'_i(L_i,X_i)&\le 0\\
		&\bar{R}'(in)&\ge\alpha,~\forall in\in|IN|\\
		&&\cdots
		\label{eq30}
	\end{array}
\end{align}
where in $\bar{R}'(in)$, the term $C_{i,k}x_{i,k}$ is replaced with $a_{i,k}L_{i,k}+b_{i,k}x_{i,k}$, hence $L_{i,k}$ and $x_{i,k}$ are also decoupled.
\subsection{Generalized Benders Decomposition}
\begin{figure}[t]
	\centerline{\includegraphics[width=0.8\linewidth]{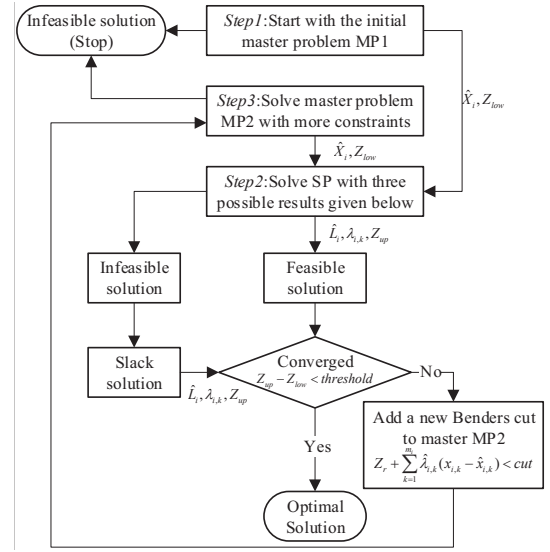}}
	\setlength{\belowcaptionskip}{-0.6cm}  
	\caption{Flowchart of Benders decomposition.}
	\label{alg}
\end{figure}

we use GBD to solve MILP problem as it can guarantee a global optimum as long as the problem is convex for fixed values of the binary variables \cite{b19}. 
The general idea of GBD is that decouple Original Problem (OP) into a Linear Programming (LP) Master Problem (MP) and a Integer Programming (IP) SubProblem (SP). By iteratively adding new cutting surface constraints, the upper and lower bounds of the objective function value are approximated, and finally converge to the optimal solution.
We take the problem \eqref{eq25} and \eqref{eq26} as an example to illustrate the algorithm, and the solution method for \eqref{eq27} and the main problem is similar. Fig.~\ref{alg} shows the proposed solution steps which are discussed below.

\textit{Step 1 (Initialization)}: Initialize the iteration counter $\ell$ = 1. We start with the initial master problem (MP)1. Solving MP1 we can get the result $\hat{X}_i^{(\ell)}$. $Z_{low}^{(\ell)}$ represents lower bound of original problem (OP).
\begin{equation}
	\begin{split}
		\min_{X_i}&\quad Z_{low}=\gamma_iB_iX^\top_i-\gamma_iN_i\\
		s.t.&\quad X_i\in \{0,1\}
		\label{eq31}
	\end{split}
\end{equation}

\textit{Step 2 (subproblem (SP) solution)}: Solve SP to get upper bound of OP. To consider possible infeasibility of the subproblem, we formulate an alternative always-feasible SP.
\begin{equation}
	\begin{split}
	\min_{L_i,\lambda_{i,k}}\quad&Z_r=L_i(\gamma_iA_iE+\rho E)+M(\sum_{k=1}^{m_i}e_k+w)\\
		s.t.\quad&L_i(\gamma_iA_iE+B_iX^\top_i)- N_i+w\ge 0\\
		&\qquad L_{i,k}-\varepsilon L^w_{i,k}x_{i,k}-e_k<0\\
		&\qquad \quad L^w_{i,k}x_{i,k}-L_{i,k}<0\\
		&\qquad \quad x_{i,k}=\hat{x}_{i,k}^{(\ell)}:\quad \lambda_{i,k}^{(\ell)}
	\label{eq32}
	\end{split}
\end{equation}
where $e_k$ and $w$ are positive slack variables to guarantee the feasibility of SP, and $M$ is a large enough positive constant. The solution of this problem is $\hat{L}^{(\ell)}_i$ with dual variable values $\hat{\lambda}_{i,1}^{(\ell)},\cdots,\hat{\lambda}_{i,m_i}^{(\ell)}$. We denote $Z_{up}^{(\ell)}$ is the upper bound of objective function value of OP.
\begin{equation}
	Z_{up}^{(\ell)}=B_i\hat{X}^{(\ell)\top}_i-N_i+Z_r^{(\ell)}
\end{equation}

If $Z^{(\ell)}_{up}-Z^{(\ell)}_{low}$ is smaller than a given tolerance, the iteration stops and the optimal results is $X^{\ell}_i$ and $L^{\ell}_i$. Otherwise, go to the next step.

\textit{Step 3 (MP2 solution)}: Update the iteration counter, $\ell \leftarrow \ell+1$. Add new benders cut constraint to filter useless solution domains.
\begin{equation}
	\begin{split}
		\min_{X_i}&\quad Z_{low}=\gamma_iB_iX^\top_i-\gamma_iN_i+cut\\
		s.t.&\quad Z_r^{(\sigma)} + \sum_{k=1}^{m_i}\hat{\lambda}^{(\sigma)}_{i,k}(x_{i,k}-\hat{x}^{(\sigma)}_{i,k})<cut\\
		&\qquad\qquad \sigma =1,\cdots,\ell-1\\
		&\qquad\qquad X_i\in \{0,1\}
	\label{eq34}
	\end{split}
\end{equation} 
Solve MP2 to obtain a new lower bound solution, then go back to \textit{Step 2} for solving SP again.

\section{Performance Evaluation}
\label{sct5}
\begin{figure*}[t]
	\begin{minipage}[t]{0.33\linewidth}
		\centerline{\includegraphics[width=1\textwidth]{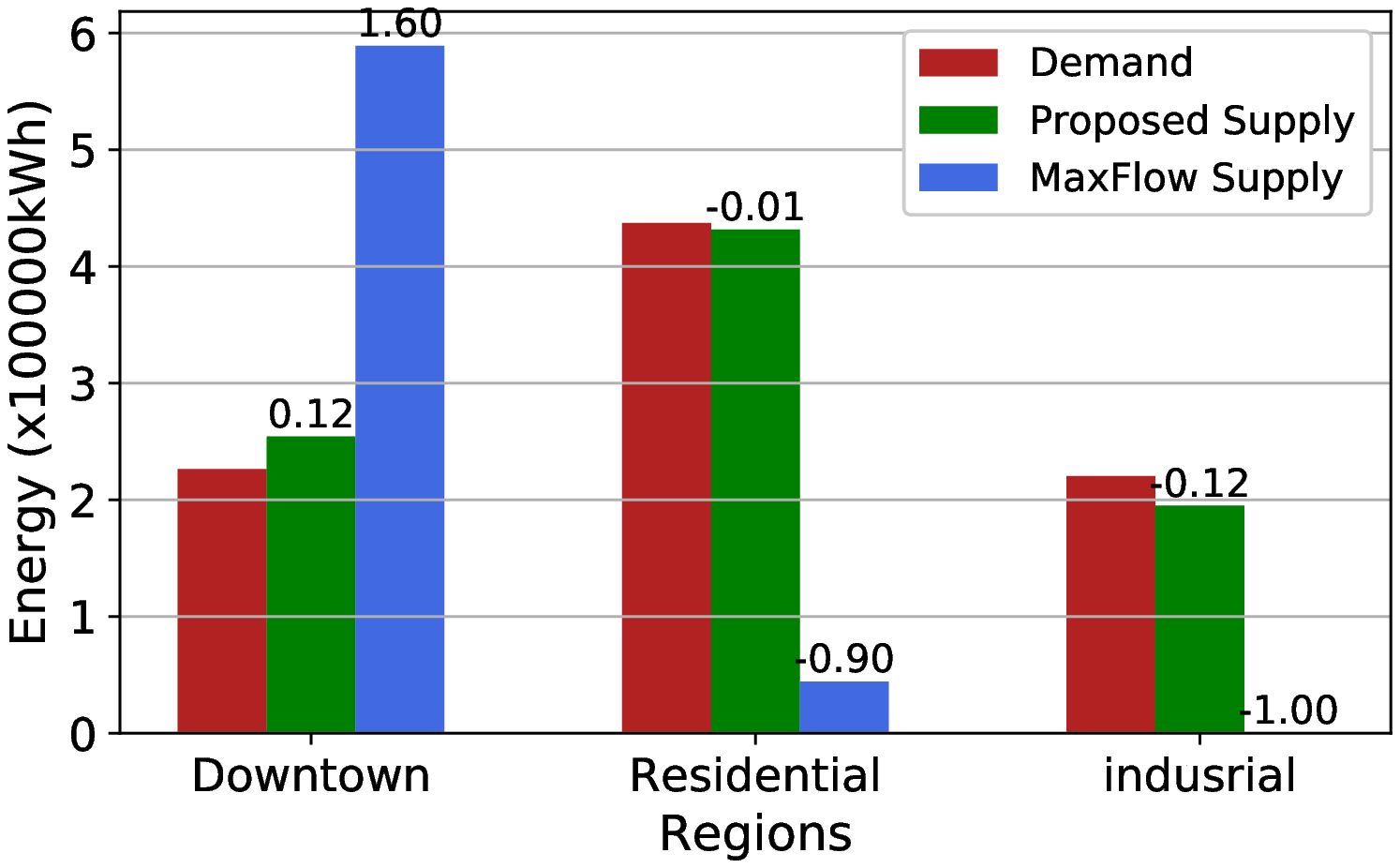}}
		\setlength{\abovecaptionskip}{-1pt}
		\caption{Energy balance.}
		\label{balance}
	\end{minipage}
	\begin{minipage}[t]{0.33\linewidth}
		\centerline{\includegraphics[width=1\textwidth]{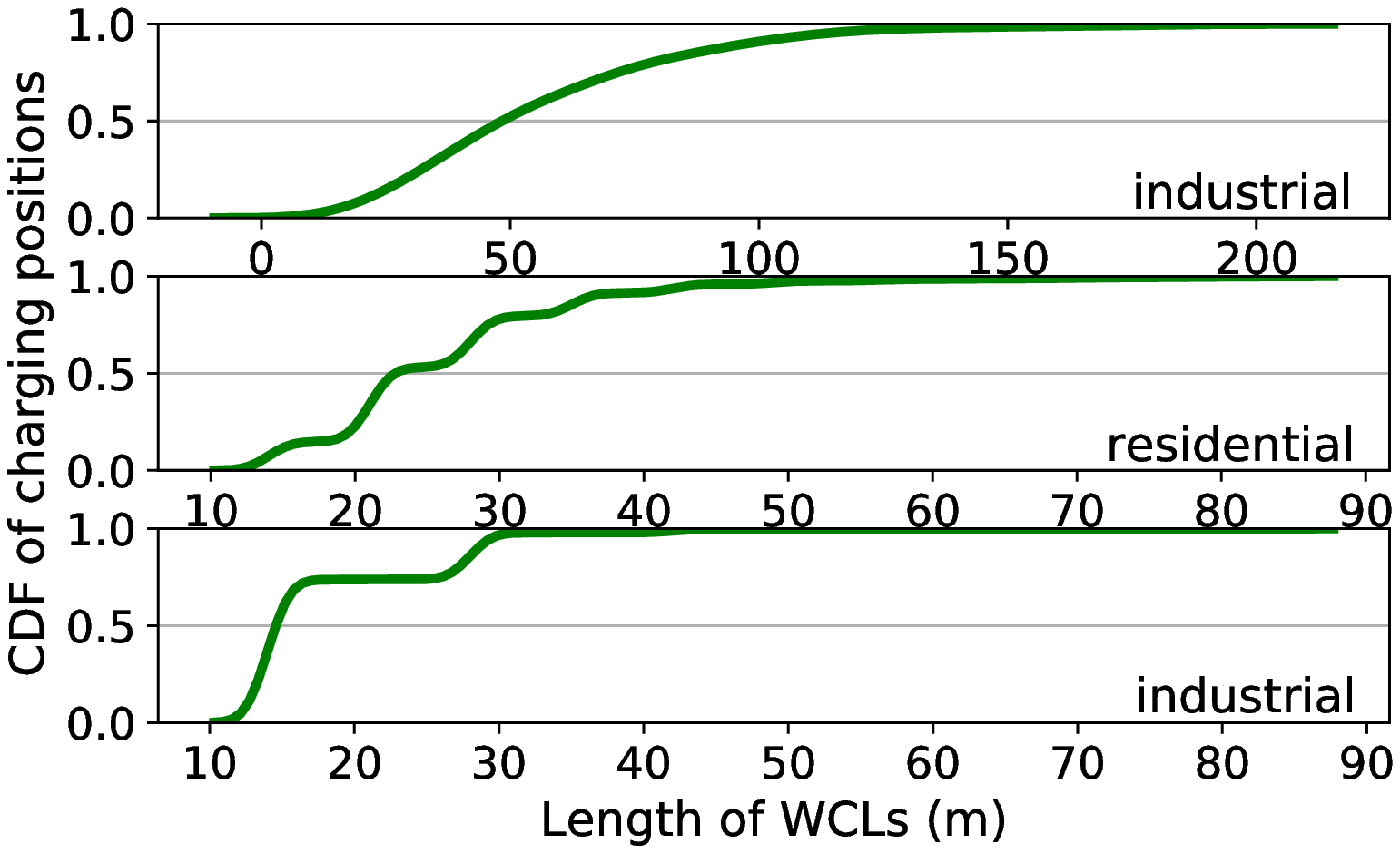}}
		\setlength{\abovecaptionskip}{-1pt}
		\caption{Distribution of WCL length.}
		\label{length_density}
	\end{minipage}
	\begin{minipage}[t]{0.33\linewidth}
		\centerline{\includegraphics[width=1\textwidth]{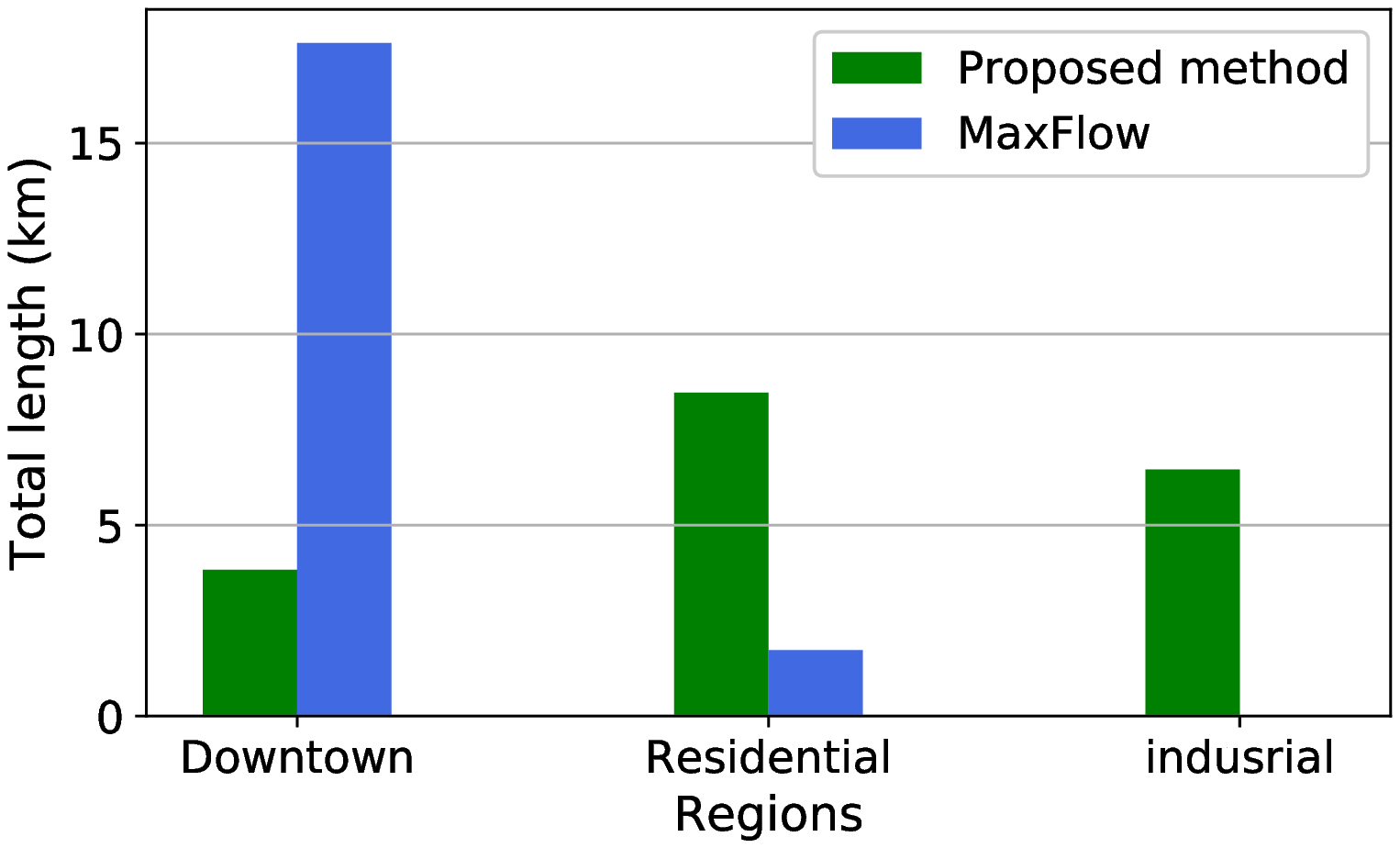}}
		\setlength{\abovecaptionskip}{-1pt}
		\caption{Total length of WCLs.}
		\label{total_length}
	\end{minipage}
	
	\begin{minipage}[t]{0.33\linewidth}
		\centerline{\includegraphics[width=1\textwidth]{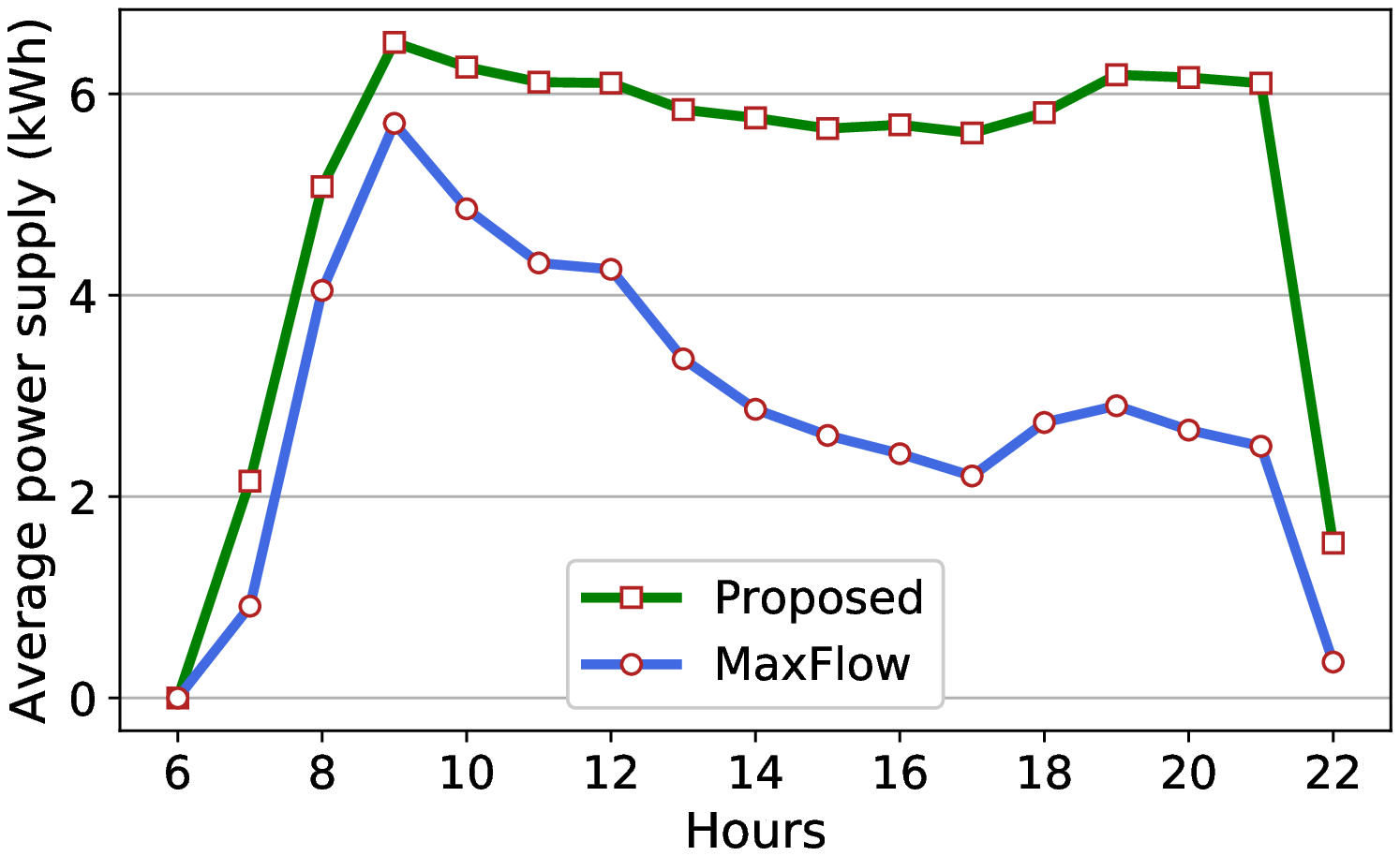}}
		\setlength{\abovecaptionskip}{-1pt}
		\setlength{\belowcaptionskip}{-1cm}   
		\caption{Average power supply.}
		\label{supply}
	\end{minipage}
	\begin{minipage}[t]{0.33\linewidth}
		\centerline{\includegraphics[width=1\textwidth]{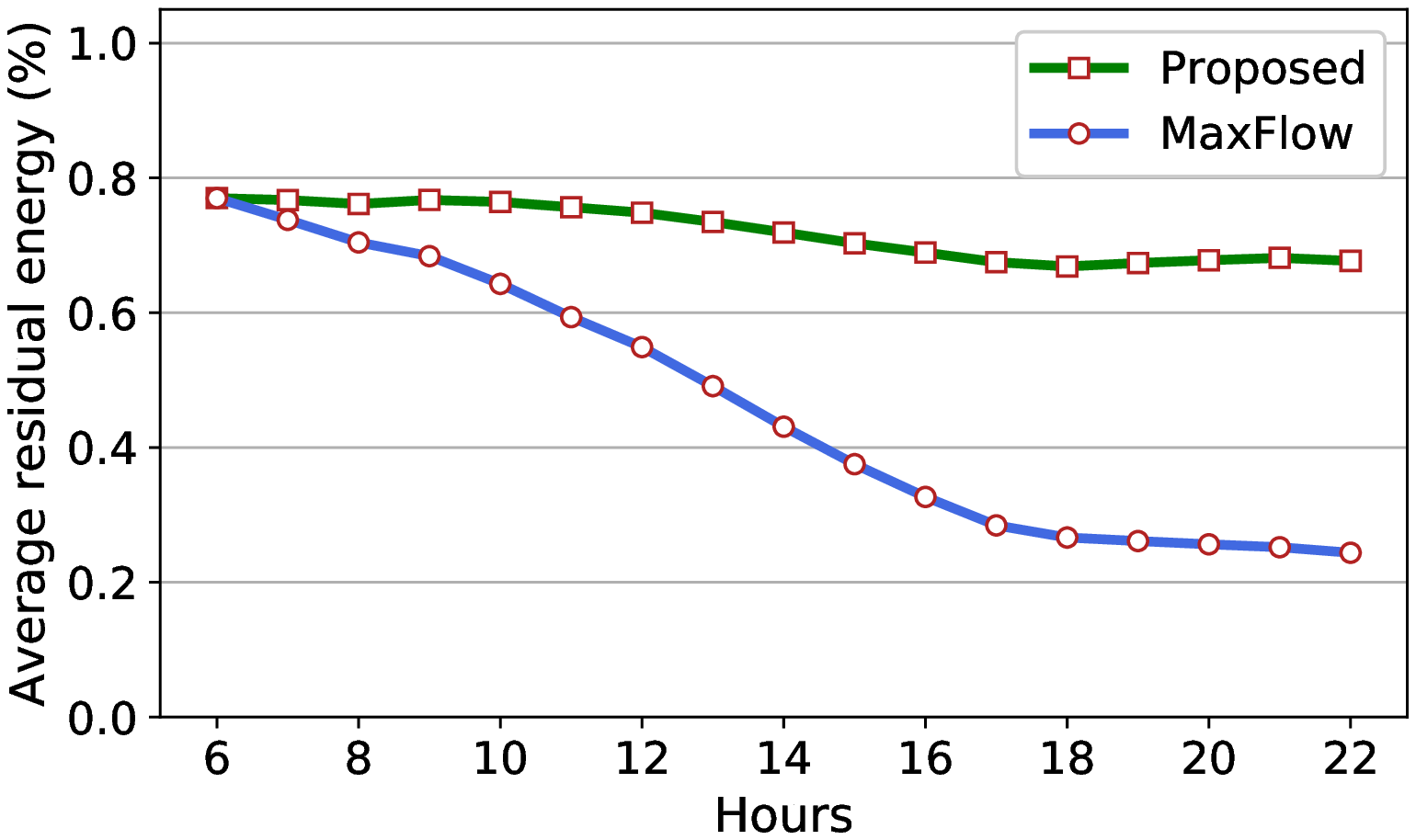}}
		\setlength{\abovecaptionskip}{-1pt}
		\setlength{\belowcaptionskip}{-1cm}  
		\caption{Average residual energy.}
		\label{residual}
	\end{minipage}
	\begin{minipage}[t]{0.33\linewidth}
		\centerline{\includegraphics[width=1\textwidth]{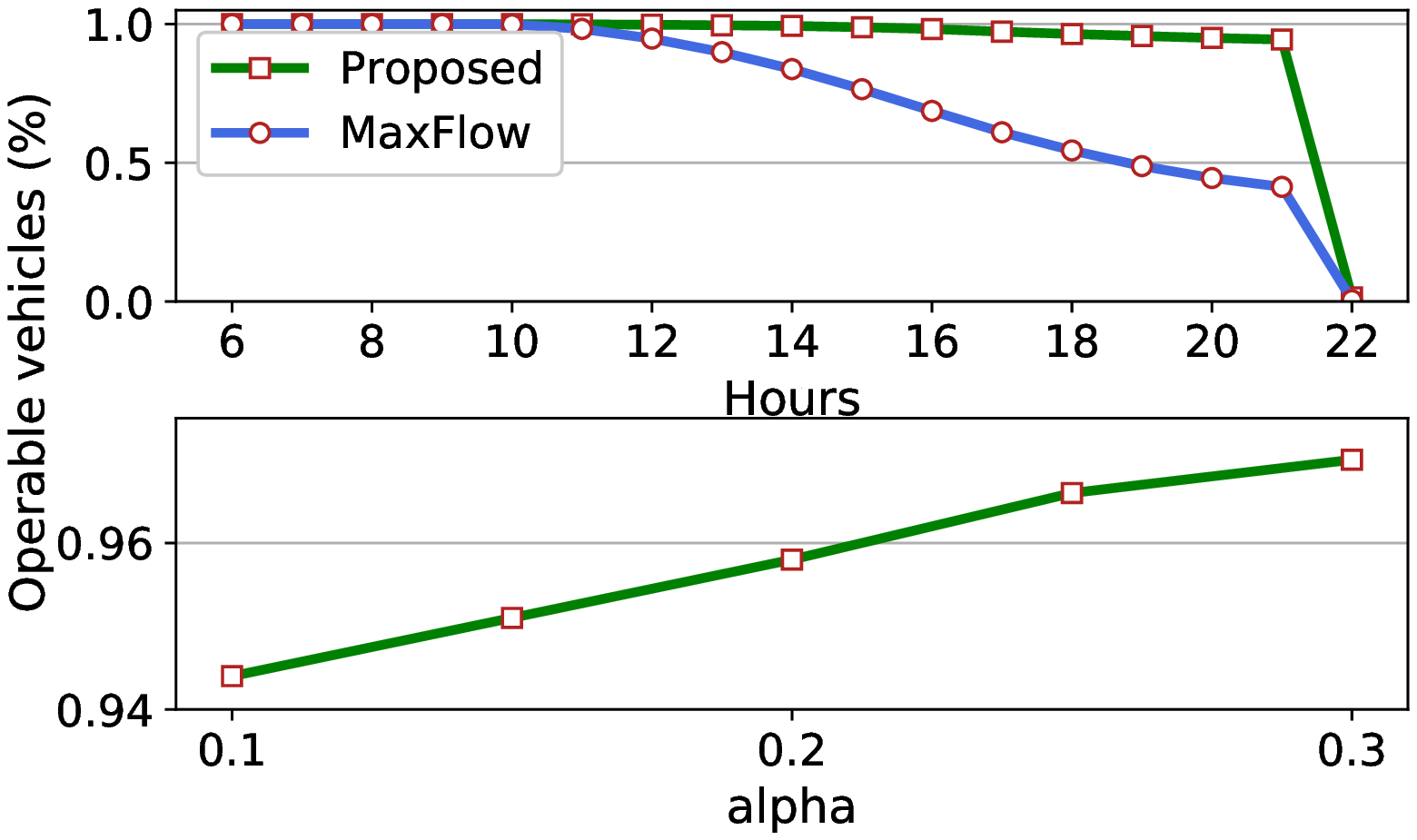}}
		\setlength{\abovecaptionskip}{-1pt}
		\setlength{\belowcaptionskip}{-1cm}  
		\caption{Operable vehicles.}
		\label{operable}
	\end{minipage}
\end{figure*}
\subsection{Experiment Settings}
First, we set up a simulation scene of urban traffic with SUMO. Since there is currently no detailed map accurate to lanes, we build a road network which covers an area of 400 $\rm km^2$ and consists of 19680 lanes and 1681 intersections. 10000 vehicles are sent into road between 6:00 to 7:00 and leave between 21:00 and 22:00. The whole urban area is divided into 3 functional subregions, and their traffic volume distributions are similar to those of the real subregions.

In simulation, the maximal SOC of vehicles are 40 kWh, and the initial SOC are randomly distributed between half and full SOC, because vehicles will be charged when the residual power of previous day is lower than half-value in case of insufficient endurance. The energy consumption rate in downtown, residential, industrial region are 0.24 kWh/km, 0.21 kWh/km, 0.19 kWh/km respectively, and the power supply rate of WCL is 100kW with 80\% power transfer efficiency\cite{b5}. The unit price of a charging lane is \$100/m \cite{b2}. Since there are no previous methods having handled the WCL deployment in urban area with TL, we created a method that maximally covers traffic flows (denoted by \textit{MaxFlow}) \cite{b6} to compare with the proposed method: the WCLs are deployed at lanes with the largest traffic flows sequentially until the deployment cost equals the one in the proposed method, and the length of each WCL is 192m (the maximum length of WCLs in proposed method). The other parameters are set as follows: $v_{equal}=$ 8m/s, $v_{min}=$ 4m/s, $\bar{d_2}=$ 7m, $\varepsilon=$ 2, $\alpha=$ 0.1, $\beta=$ 1, $\gamma_1=10^4$, $\gamma_2=10^5$, $\gamma_3=10^4$, $\eta\%=30\%$.
\subsection{Experiment Results}
\subsubsection{Regional Energy Balance}
To verify the correctness of our power supply model and the effect of regional energy supply and demand balance, energy supply are compared with charging demand in each subregion. Fig.~\ref{balance} shows that the total supply of proposed method and the total demand in the whole urban areas are basically the same. The balance errors in downtown and industrial region are less than 12\%, and the one in residential region is less than 1\%, which conforms to the feature of demand fluctuation. But the MaxFlow method cannot achieve high-quality energy supply, since lanes with high traffic flow are mainly concentrated in downtown region. Under the same limitation of total installation length, MaxFlow cannot meet the energy demand of other subregions.
\subsubsection{Length of WCLs}
Fig.~\ref{length_density} indicates the Cumulative Distribution Function (CDF) of WCLs` length in each subregion of proposed method. The MaxFlow is not included because length of almost all WCLs are the same. In downtown region, more than 50\% WCLs are longer than 50m and the curve is smooth, for the reason that traffic condition in this subregion is generally poor. The length distribution of residential areas is clearly multi-stratified, because the traffic conditions in this subregion are quite rich. While in industrial region, more than 75\% WCLs are shorter than 20m and the number of WLCs is relatively large, since roads in this subregion are relatively smooth, and vehicles do not queue at intersections. More WCLs are needed to supply enough energy to EVs.

Fig.~\ref{total_length} shows the comparison of the total length of WCLs in each subregion. WCLs of MaxFlow are concentrated in downtown region as mentioned above. Comparing Fig.~\ref{balance} and Fig.~\ref{total_length} we can find that, in downtown region, total WCL lengths of MaxFlow beyond the one of proposed method cannot lead to an equal proportion of energy increase, owing to unreasonable distribution of MaxFlow. While in industrial subregion, the energy supply per unit length of WCLs is significantly less than that of the other two subregions, due to the smooth traffic condition.  
\subsubsection{Energy Condition of Vehicles}
As Fig.~\ref{supply} indicates, since shared EVs start during 6:00-7:00, and gradually stop working between 21:00 and 22:00, the hourly average power supply for each vehicle in both durations are relatively low. Increasing vehicles are going to downtown in rush hours 7:00-9:00, so the power supply of both methods increase rapidly. But during 9:00-17:00, EVs tend to flow to other subregion to seek passengers. Since WCLs of MaxFlow are concentrated in downtown region, its supply declines rapidly while the one of proposed method drops slowly. During 17:00-20:00, downtown and residential regions attract some operable vehicles, so the supply is a bit increased.

Fig.~\ref{residual} shows the average residual energy of vehicles under different hours. It is clear that the residual energy of proposed method is larger than the one of MaxFlow due to the same reason of average power supply condition.
\subsubsection{Operable Vehicles}
Fig.~\ref{operable} shows that at the beginning of 6:00-10:00, the number of operable vehicles decreases slowly due to the initial battery quantity is high. After 10:00, the number of operable vehicles of MaxFlow drops faster than that of proposed method, for the reason that many EVs of MaxFlow cannot enjoy dynamic charging service in subregions except downtown. While in proposed method, most EVs are able to find WCLs on their inter-region trips since we have considered the recharge-consumption ratio from global view.

The variation of the number of drivable vehicles with parameters $\alpha$ at the end of a day is also shown in Fig.~\ref{operable}. With the increase of $\alpha$ the operable vehicles increase linearly,  because greater $\alpha$ means that vehicles are more likely to get recharged in time when they shuttle back and forth between distant WCLs. Actually these simulation is without considering recharge-consumption ratio of industrial region, otherwise it is infeasible, because power supply capacity of WCLs in industrial are too low to satisfy this constraint. 
\section{Conclusion}
\label{sct6}
In this paper, we investigate the characteristics of urban charging demand based on real metropolitan traffic datasets, finding that the average daily energy demand of taxis in the whole urban areas is basically unchanged, while regional demand shows varying degrees of fluctuation. To obtain deployment scheme conforming to these features, we build a WCL power supply model, which is closely related to traffic conditions including TLs, to characterize the power supply of the WCLs. Bilevel utility functions are formulated to achieve the lowest construction cost and balance of supply and demand in each subregion as far as possible considering EV operability, under the premise of matching supply and demand in the whole urban areas. Comparing with the simulation results of the baseline, we found that placing all WCLs at intersections with large traffic volume is not the best choice. The proposed method can better utilize WCL construction cost and achieve balance of supply and demand in downtown and residential regions. In addition, we also suggest that industrial region is not suitable for deploying WCLs due to its high cost-supply ratio and low satisfaction to vehicle operability. 
\section*{Acknowledgment}
The authors are of the Department of Automation, Shanghai Jiaotong University, and Key Laboratory of System Control and Signal Processing, Ministry of Education of China, Shanghai 200240, China. This work was supported by National Key Research and Development Program of China (2016YFB0901900) and National Natural Science Foundation of China (61573245, 61731012, 61633017, and 61622307).
\bibliographystyle{IEEEtran}
\bibliography{ref}
\end{document}